\documentclass[12pt]{amsart}
\usepackage{amsmath}
\usepackage{amssymb}
\usepackage{amsmath,amssymb,amsthm,amscd}

\usepackage{txfonts}
\usepackage{amsbsy}
\usepackage{graphicx,color}

\numberwithin{equation}{section}

\newtheorem{prop}{Proposition}[section]
\newtheorem{theo}{Theorem}[section]
\newtheorem{lemm}{Lemma}[section]
\newtheorem{coro}{Corollary}[section]

\newtheorem{defi}{Definition}[section]

\def\begeq{\begin{equation}}
\def\endeq{\end{equation}}

\begin{document}

\title{Energy non-collapsing and refined blowup for a semilinear heat equation}
\author{Shi-Zhong Du}
\thanks{$^\dagger $The paper is partially supported by STU Scientific Research Foundation for Talents (SRFT-NTF16006).}
\address{The Department of Mathematics,
            Shantou Univeristy, Shantou, 515063, P. R. China.} \email{szdu@stu.edu.cn}

\renewcommand{\subjclassname}{%
  \textup{2000} Mathematics Subject Classification}
\date{Apr. 2018}

\begin{abstract}
   Refined structures of blowup for non-collapsing maximal solution to a semilinear
    parabolic equation
       $$
        u_t-\triangle u=|u|^{p-1}u
       $$
with $p>1$ are studied. We will prove that the blowup set is empty for non-collapsing blowing-up in subcritical case, and all finite time non-collapsing blowing-up must be type II in critical case. When $p>p_S\equiv\frac{N+2}{N-2}$ for $N\geq3$, the Hausdorff dimension of the blowup set for maximal solution whose energy is non-collapsing is shown to be no greater than $N-2-\frac{4}{p-1}$, which answers a question proposed in \cite{CDZ} positively. At the end of this paper, we also present some new examples of collapsing and non-collapsing blowups.
\end{abstract}

\maketitle\markboth{Shi-Zhong Du}{Energy Non-collapsing}

\section{Introduction}

 Consider the following classical semilinear parabolic equation
 \begin{equation}\label{e1.1}
        u_t-\triangle u=|u|^{p-1}u \ \mbox{ in }
        \Omega\times[0,\omega)
 \end{equation}
under Dirichlet boundary condition
    $$
      u=0\ \ \ \mbox{ on } \partial\Omega\times(0,\omega),
    $$
where $\Omega$ is a bounded domain in ${\mathbb{R}}^N$ with smooth
boundary and the nonlinear exponent $p$ is assumed to be greater than one. It's well known that for any initial datum in $L^\infty(\Omega)$, this problem admits a unique classical solution in
${\Omega}\times(0,\omega)$ for $0<\omega\leq\infty$, which
is maximal in the sense that
     $$
       \limsup_{t\uparrow\omega}||u(t)||_{L^\infty(\Omega)}=\infty
     $$
when $\omega$ is finite. After the pioneering works of Kaplan-Fujita-Levine-Ball et al \cite{K,F,L,B}, one knows blowup solutions with some special initial data. For example, when the initial energy is negative, the solution blows up in finite time. In order to understand the blowup phenomenons, we need to explore further on different aspects. One of these topics is the completeness of blowup which was proposed by Brezis \cite{BC}, concerning the extendability of maximal solutions who blow up in finite beyond the first singular time. A first result was obtained by Baras-Cohen in \cite{BC}, where they showed that the blowup must be complete under either one of the following conditions:

(A1) $\Omega$ is convex and $u_0\geq0, \triangle u_0+u_0^p\geq0$, or

(A2) $\Omega$ is convex, $u_0\geq0$ and $1<p<p_S$.

Later, Galaktionov-Vazquez \cite{GV} extended the result of \cite{BC} to supercritical case for some special radial symmetric solutions.

 Another striking question concerns the blowup rate of the singularity. Like that in curvature flow, it's natural to distinguish the singularities into two classes. The first one is the so-called Type I blowup in the sense
    \begin{equation}\label{e1.2}
       \limsup_{t\to
       \omega^-}\sup_{x\in\Omega}(\omega-t)^{\frac{1}{p-1}}|u(x,t)|<\infty,
    \end{equation}
and another one is the Type II blowup, which means that
    \begin{equation}\label{e1.3}
       \limsup_{t\to
       \omega^-}\sup_{x\in\Omega}(\omega-t)^{\frac{1}{p-1}}|u(x,t)|=\infty.
    \end{equation}
Type I blowing-up rate which comes from the trivial solution
     $$
       u(t)=\kappa(\omega-t)^{-\frac{1}{p-1}}
     $$
 of\eqref{e1.1} with $\kappa=(p-1)^{-\frac{1}{p-1}}$, is sometimes characterized by the self-similar blowing-up \cite{GK1}. Weissler proved in \cite{W2} that finite time blowup must be type I under radial symmetric assumption and some rather special initial data. At the same time, Friedman-McLeod \cite{FM} have also derived a similar result for convex $\Omega$ and monotone solution $u(\cdot,t)$ in $t$. Recently, Matano-Merle \cite{MM} obtained a satisfactory result under radial symmetric case, which states that for any $p$ lies between $p_S$ and the Joseph-Lundgren \cite{JL} exponent
   $$
    p_{JL}=\begin{cases}
      \infty, & \mbox{ if } 3\leq N\leq10,\\
      1+\frac{4}{N-4-2\sqrt{N-1}}, & \mbox{ if } N\geq11,
    \end{cases}
   $$
 all radial symmetric solutions can only blow up in finite time with type I rates, and for $p=p_S$, a similar conclusion holds for all positive radial symmetric solutions. The upper bound $p_{JL}$ of $p$ in his theorem can not be improved due to the existence of type II blowup in region of $p>p_{JL}$ \cite{HV1,HV2}. Result of type I blowup rate without radial symmetric was firstly discovered in \cite{GK1}, where Giga-Kohn showed that when $\Omega$ is a convex domain or whole space ${\mathbb{R}}^N$, finite time blowup must be type I under one of the following assumptions:

 (A3) $u_0\geq0$ and $1<p<p_S$, or

 (A4) $u_0$ may change sign and $1<p<\frac{3N+8}{3N-4}$ for $N\geq2$ or $p>1$ for $N=1$.

\noindent This result was later extended by Giga-Matsui-Sasayama \cite{GMS} to all subcritical exponents and all change sign solutions for Cauchy problem.

In this paper, we will prove the following result concerning the first singular time $\omega<+\infty$ under critical and subcritical case:\\

 \noindent\textbf{Theorem A.} Assume that $1<p<p_S$, no non-collapsing blowup can occur. Hence, all finite time blowup must be complete under subcritical case. If $p=p_S$, all non-collapsing blowup must be type II in the sense \eqref{e1.3}.\\

 Unlike that in subcritical case, after work of \cite{CDZ}, one knows that positive borderline solutions are all examples of non-collapsing blowup in supercritical case (see also \cite{MM2} for some other examples of non-collapsing blowup). Furthermore, it was proved in \cite{CDZ} that for convex domain and all supercritical exponents $p$, the singular set of borderline solution is not empty and of finite $\Big(N-\frac{4}{p-1}\Big)-$Hausdorff measure. The convex assumption was later weaken to star-shaped domain in \cite{CD}. Concerning the Hausdorff dimension of the singular set at each time slice, the authors also laid down a conjecture in \cite{CDZ} by claiming that it may be possible to lower the Hausdorff dimension to $N-\frac{4}{p-1}-2$. In this paper, we will give a positive answer to it in term of non-collapsing blowup:\\

\noindent\textbf{Theorem B.}   Assume that $p>p_S$ and $\Omega$ is convex. Let $u$ be a maximal solution which blows up non-collapsing at $t=\omega$ and ${\mathcal{S}}$ be the compact blowup set whose Hausdorff dimension is no greater than $N-\frac{4}{p-1}-2$. We have $u$ is locally bounded near $(a,\omega)$ when $a$ doesn't belong to ${\mathcal{S}}$.\\

It's worthy to compare our result with that in \cite{V}, where Velazquez has shown that the Hausdorff dimension of the blowup set for positive solution to Cauchy problem of \eqref{e1.1} with subcritical nonlinearity must be less than or equal to $N-1$. (A similar result was also obtained in \cite{MM} for all radial symmetric solutions when $p_S<p<p_{JL}$ and for all positive radial symmetric solutions when $p=p_S$.) Recently, Blatt-Struwe independently derived in \cite{BS} (see Theorem 3.8 there) that the Hausdorff dimension of the blowup set can not exceed $N-\frac{4}{p-1}-\kappa$ for some $1\leq\kappa\leq2$ under an extra assumption
   \begin{equation}\label{e1.5}
      \liminf_{R\downarrow 0}R^{-\kappa}\int^\omega_{\omega-R^2}\int_\Omega(|\nabla u|^2+|u|^{p+1})dxdt<\infty.
   \end{equation}
Combining with our Theorem \ref{t2.1} below, it's not hard to see that for non-collapsing blowup, \eqref{e1.5} holds for $\kappa$ closing to $2$ arbitrarily from below, and thus gives another proof to our Theorem B.

 The third part of the paper is devoted to construction of new examples of complete and incomplete blowups. A first complete blowup solution was found by Baras-Cohen in \cite{BC}, where they showed that for $1<p<\frac{N}{N-2}$ when $N\geq3$ or for $p>1$ when $N=2$, the maximal solution blows up completely with some rather special initial data $u_0$. Later, the result was extended to all subcritical and critical exponents $1<p\leq p_S$ for radial symmetric solutions by Galaktionov-Vazquez in \cite{GV}. In this paper, we shall prove the following results of complete and incomplete blowups by applying our refined criteria for non-collapsing singularities:\\

 \noindent\textbf{Theorem C.}  For $1<p<p_S$, all finite time blowup must be complete.

 If $p=p_S$ and $\Omega$ is convex, type I blowup must be collapsing. Therefore, all finite time blowup must be complete in case $\Omega$ is a ball and $u_0$ is a radial symmetric positive function.

 When $p>p_S$, the blowup must be non-collapsing under either one of the following conditions:

(A5)  $u$ is a borderline solution to \eqref{e1.1}, whose finite time blowup was guaranteed in \cite{CD,CDZ} for star-shaped domain $\Omega$, or

(A6) $\Omega=B_R, R>0$, $u_0=u_0(r)>0$ and $u'_0(r)\leq0, \forall 0<r<R$.

On the other hand, the blowup must be collapsing and hence complete under the following condition:

(A7) $\Omega=B_{R}\setminus B_{R_0}$, $u_0=u_0(r)$ for all $R_0<r=|x|<R$.\\

Our full complete blowup result for subcritical case improves that of \cite{BC} and \cite{GV}. Under critical case, although our result is a special case of \cite{GV}, we give a different proof to complete blowup without Zero Comparison Theorem, which has a potential to construct examples of complete blowup that is not radial symmetry. Finally, some new complete and incomplete blowups were observed for supercritical case.

\vspace{40pt}

\section{Energy collapsing and complete blow-up}

Let $u(\cdot,t)$ be a maximal solution to \eqref{e1.1} on $[0,\omega)$. Multiplying \eqref{e1.1} by $u_t$ and integrating over $\Omega$, we get
  \begin{equation}\label{e2.1}
    \frac{d}{dt}E(u)=-\int_\Omega u_t^2, \ \ \ \forall t\in[0,\omega)
  \end{equation}
after integration by parts, where
  $$
    E(u)\equiv\frac{1}{2}\int_\Omega|\nabla u|^2-\frac{1}{p+1}\int_\Omega|u|^{p+1}
  $$
is called to be the energy of the solution. Since the energy is monotone non-increasing by \eqref{e2.1}, it's natural to define its limit by
   \begin{equation}\label{e2.2}
     {\mathcal{B}}\equiv\lim_{t\to\omega^-}E(u(t)),
   \end{equation}
which is finite or negative infinity.

\begin{defi}\label{d2.1}
  We call the energy (or blowup) to be collapsing if ${\mathcal{B}}=-\infty$ and to be non-collapsing  if ${\mathcal{B}}>-\infty$.
\end{defi}

Also, we can define the complete blow-up in spirit of \cite{GV}:
\begin{defi}\label{d2.2}
  Let $u$ be a maximal solution to \eqref{e1.1} which blows up at time $t=\omega$. Suppose that there exist $\delta>0$ and a sequence of classical solutions $\{u_k\}_{k=1}^\infty$ of \eqref{e1.1} defined on $\overline{\Omega}\times[0,\omega+\delta)$, which tends to $u$ uniformly on any compact subset of $\overline{\Omega}\times[0,\omega)$, we will call the blowup to be incomplete.

  Otherwise, we call it to be complete blowup.
\end{defi}

To clarify the concept of incomplete blow-up more clearly, we have the following Proposition:

\begin{prop}\label{p2.1}
  Let $p>1$ and $u$ be a maximal solution to \eqref{e1.1} which blows up incompletely at $t=\omega$. Then there exist a $\delta>0$ and a weak solution $U$ of \eqref{e1.1} on $\Omega\times[0,\omega+\delta/2)$, such that $U$ is identical to $u$ on $\Omega\times[0,\omega)$.
\end{prop}

\noindent\textbf{Proof.} Let $\delta>0$ and $\{u_k\}_{k=1}^\infty$ be given in Definition \ref{d2.2}. To show the existence of weak solution on $\Omega\times[0,\omega+\delta/2)$, we need only to prove some uniform a-priori bounds on $u_k$. Setting $v=u_k$ for simplicity, multiplying \eqref{e1.1} by $v$ and then integrating over $\Omega$, we get
   \begin{equation}\label{e2.3}
     \frac{1}{2}\frac{d}{dt}\int_\Omega v^2=-\int_\Omega|\nabla v|^2+\int_\Omega|v|^{p+1}=-2E(v)+\frac{p-1}{p+1}\int_\Omega|v|^{p+1}
   \end{equation}
after integration by parts. Claim: There exists a positive constant $C_{p,\Omega}$ depending only on $p$ and $\Omega$, such that
  \begin{equation}\label{e2.4}
     E(v_0)\geq E(v(t))\geq-C_{p,\Omega}(\omega+\delta-t)^{-\frac{p+1}{p-1}}
  \end{equation}
for any $0<t<\omega+\delta$. Without loss of generality, we may assume that $E(v(t_0))<0$ for some $0<t_0<\omega+\delta$. Then $E(v(t))<0$ for all $t_0<t<\omega+\delta$. Noting first that
   \begin{eqnarray}\label{e2.5}\nonumber
     \frac{p-1}{p+1}\int_\Omega|v|^{p+1}+2|E(v)|&=&\int_\Omega vv_t\leq\Bigg(\int_\Omega|v|^{p+1}\Bigg)^{\frac{1}{p+1}}\Bigg(\int_\Omega |v_t|^{\frac{p+1}{p}}\Bigg)^{\frac{p}{p+1}}\\
      &\leq&\frac{p-1}{2(p+1)}\int_\Omega|v|^{p+1}+C_{p,\Omega}\Bigg(\int_\Omega v_t^{2}\Bigg)^{\frac{p+1}{2p}}, \ \ \forall t>t_0
   \end{eqnarray}
by \eqref{e2.3} together with H\"{o}lder and Young's inequalities, there exists some positive constant $\delta_{p,\Omega}$ depending only on $p$ and $\Omega$, such that
   \begin{equation}\label{e2.6}
      \int_\Omega v_t^2\geq \delta_{p,\Omega}|E(v)|^{\frac{2p}{p+1}}
   \end{equation}
for all $t_0<t<\omega+\delta$. Substituting \eqref{e2.6} into \eqref{e2.1}, we get
  \begin{equation}\label{e2.7}
    \frac{d}{dt}|E(v)|\geq \delta_{p,\Omega}|E(v)|^{\frac{2p}{p+1}}, \ \forall t\in[t_0,\omega+\delta).
  \end{equation}
Therefore, \eqref{e2.4} follows by solving the ordinary differential inequality \eqref{e2.7} in time. Now, integrating \eqref{e2.1} over $[0,\omega+\delta/2]$ and using \eqref{e2.4}, we conclude that
   \begin{equation}\label{e2.8}
     \int^{\omega+\delta/2}_0\int_\Omega v_t^2\leq C(p,\Omega,\delta,\omega,E(v_0)).
   \end{equation}
Claim: There exists a positive constant $C'_{p,\Omega}$ such that
  \begin{equation}\label{e2.9}
     \int_\Omega v^2(t)\leq C'_{p,\Omega}\max\Big\{|E(v_0)|^{\frac{2}{p+1}}, (\omega+\delta-t)^{-\frac{2}{p-1}}\Big\}
  \end{equation}
for any $0<t<\omega+\delta$. In fact, assuming that for some $0<t_0<\omega+\delta$,
   $$
     \int_\Omega v^2(t_0)>\Bigg\{\frac{4(p+1)}{p-1}|\Omega||E(v_0)|\Bigg\}^{\frac{2}{p+1}},
   $$
we derive from \eqref{e2.3} and H\"{o}lder's inequality that
  \begin{equation}\label{e2.10}
    \frac{1}{2}\frac{d}{dt}\int_\Omega v^2\geq\frac{p-1}{2(p+1)}|\Omega|^{-\frac{p-1}{2}}\Bigg(\int_\Omega v^2\Bigg)^{\frac{p+1}{2}}.
  \end{equation}
Integrating over $(t_0,\omega+\delta)$, we conclude that
   $$
     \int_\Omega v^2(t)\leq C'_{p,\Omega}(\omega+\delta-t)^{-\frac{2}{p-1}}
   $$
for all $t\in[t_0,\omega+\delta)$. So, the claim is true as long as $C'_{p,\Omega}$ is taken to be greater than $\Big\{\frac{4(p+1)}{p-1}|\Omega|\Big\}^{\frac{2}{p+1}}$. Finally, we show that
   \begin{equation}\label{e2.11}
    \int^{\omega+\delta/2}_0\int_\Omega\Big(|\nabla v|^2+|v|^{p+1}\Big)dxdt\leq C(p,\Omega,\delta,\omega,E(v_0))<\infty
   \end{equation}
holds for some positive constant $C(p,\Omega,\delta,\omega, E(v_0))$ depending only on $p,\Omega,\delta,\omega, E(v_0)$. In fact, integrating \eqref{e2.3} and using \eqref{e2.4} together with \eqref{e2.9}, we have
  \begin{eqnarray*}
    &&\int^{\omega+\delta/2}_0\int_\Omega\Big(|\nabla v|^2+|v|^{p+1}\Big)dxdt\\
     &&\ \ \ \ \ \ \ \ \ \ \ \ \ \ \ \  \ \ \ \ \ \ \ \leq 2E(v_0)(\omega+\delta/2)+\frac{p+3}{p+1}\int^{\omega+\delta/2}_0\int_\Omega|v|^{p+1}dxdt\\
    &&\ \ \ \ \ \ \ \ \ \ \ \ \ \ \ \  \ \ \ \ \ \ \ \leq\frac{4(p+1)}{p-1}E(v_0)(\omega+\delta/2)+\frac{p+3}{2(p-1)}\int_\Omega v^2(\omega+\delta/2)\\
    &&\ \ \ \ \ \ \ \ \ \ \ \ \ \ \ \  \ \ \ \ \ \ \  \leq C(p,\Omega,\delta,\omega,E(v_0))<\infty.
  \end{eqnarray*}
Now, A combination of \eqref{e2.8} with \eqref{e2.9}, \eqref{e2.11} gives the desired uniform a-priori estimations for $v=u_k$, and thus a limiting weak solution $U$ on $\Omega\times[0,\omega+\delta/2)$ in the sense of distribution by passing to the limits. $\Box$\\

We also have the following relationship between energy collapsing and complete blowup:

\begin{prop}\label{p2.2}
  For all $p>1$, if the energy is collapsing, then the blow-up is complete.
\end{prop}

It would be interesting to ask whether or not the blowup is incomplete when the energy is non-collapsing. To prove the proposition, we need the following Lemma:

\begin{lemm}\label{l2.1}
  Let $v$ be a classical solution to \eqref{e1.1}. Suppose that the energy $E(v(t_0))$ becomes negative for some $t_0>0$, then the solution must blow up before
  $$
    t_0+C''_{p,\Omega}\big|E(v(t_0))\big|^{-(p-1)/(p+1)},
  $$
where $C''_{p,\Omega}$ is a positive constant depending only on $p$ and $\Omega$.
\end{lemm}

\noindent\textbf{Proof.} Noting first that the energy is monotone non-increasing, so $E(v(t))\leq E(v(t_0))<0$ for all $t>t_0$. Multiplying \eqref{e1.1} by $v$ and integrating over $\Omega$, we get
  \begin{eqnarray}\label{e2.12}\nonumber
    \frac{1}{2}\frac{d}{dt}\int_\Omega v^2&=&-2E(v)+\frac{p-1}{p+1}\int_\Omega|v|^{p+1}\\ \nonumber
      &\geq&2\big|E(v(t_0))\big|+\frac{p-1}{p+1}|\Omega|^{-\frac{p-1}{2}}\Bigg(\int_\Omega v^2\Bigg)^{\frac{p+1}{2}}\\
      &\geq& \delta_{p,\Omega}\Bigg(\big|E(v(t_0))\big|^{\frac{2}{p+1}}+\int_\Omega v^2\Bigg)^{\frac{p+1}{2}}
  \end{eqnarray}
after integration by parts, where $\delta_{p,\Omega}$ is some positive constant depending only on $p$ and $\Omega$. Setting
   $$
     y(t)\equiv\big|E(v(t_0))\big|^{\frac{2}{p+1}}+\int_\Omega u^2(t)
   $$
and rewriting \eqref{e2.12} by
   $$
     \frac{dy^{\frac{1-p}{2}}}{dt}\leq-(p-1)\delta_{p,\Omega},
   $$
we get
   $$
     0<\Bigg(\big|E(v(t_0))\big|^{\frac{2}{p+1}}+\int_\Omega v^2(t_0)\Bigg)^{\frac{1-p}{2}}-(p-1)\delta_{p,\Omega}(t-t_0)
   $$
or
   $$
     t-t_0\leq(p-1)^{-1}\delta_{p,\Omega}^{-1}\Bigg(\big|E(v(t_0))\big|^{\frac{2}{p+1}}+\int_\Omega v^2(t_0)\Bigg)^{-\frac{p-1}{2}}
   $$
as long as no blowing-up occurs before time $t$. So, the solution does blow up before
   $$
     t_0+C''_{p,\Omega}\big|E(v(t_0))\big|^{-(p-1)/(p+1)}
   $$
for $C''_{p,\Omega}=2(p-1)^{-1}\delta^{-1}_{p,\Omega}$. The conclusion is drawn. $\Box$\\

 Now, Proposition \ref{p2.2} follows from Lemma \ref{l2.1}. In fact, suppose on the contrary, then there exist a $\delta>0$ and a sequence of classical solution $\{u_k\}_{k=1}^\infty$ of \eqref{e1.1} defined on $\overline{\Omega}\times[0,\omega+\delta)$, such that $u_k$ tends to $u$ uniformly on any compact set of $\overline{\Omega}\times[0,\omega)$. Let's take $t_0<\omega$ closing to $\omega$ sufficiently, such that
    $$
      E(u(t_0))<-2\Bigg(\frac{2C''_{p,\Omega}}{\delta}\Bigg)^{\frac{p+1}{p-1}}.
    $$
 By uniform convergence, we have
     $$
      E(u_k(t_0))<-\Bigg(\frac{2C''_{p,\Omega}}{\delta}\Bigg)^{\frac{p+1}{p-1}}
    $$
 for $k$ large. Therefore, applying Lemma \ref{l2.1} to $v=u_k$, we know that $u_k$ must blow up before $t_0+\frac{\delta}{2}<\omega+\delta$, which contradicts with our assumption. \\

 For all nonlinear exponent $p>1$, we have the following a-priori estimation concerning non-collapsing blowup:

 \begin{theo}\label{t2.1}
   Let $p>1$ and $u(\cdot,t)$ be a maximal solution to \eqref{e1.1} which blows up non-collapsing at time $t=\omega$. We have for any $q\geq1$, there holds
      $$
        \int^\omega_0\Bigg(\int_\Omega(|\nabla u|^2+|u|^{p+1})dx\Bigg)^qdt\leq C_q<+\infty.
      $$
  Specially, if $1<p<p_S$, we have
      $$
       \sup_{t\in[0,\omega)}\int_\Omega|u|^{p+1}(t)<\infty
      $$
  and hence the boundedness of solution up to $t=\omega$.
 \end{theo}

 \noindent\textbf{Proof.} Noting first that it follows from \eqref{e2.1} and ${\mathcal{B}}>-\infty$ that
    \begin{equation}\label{e2.13}
      \int^\omega_0\int_\Omega u_t^2dxdt\leq E(u_0)-{\mathcal{B}}.
    \end{equation}
 Multiplying \eqref{e1.1} by $u$ and integrating over $\Omega$, we get
   \begin{equation}\label{e2.14}
     \frac{1}{2}\frac{d}{dt}\int_\Omega u^2=-2E(u)+\frac{p-1}{p+1}\int_\Omega|u|^{p+1}
   \end{equation}
 after integration by parts similar to \eqref{e2.3}. Consequently, it yields from \eqref{e2.14} that
   \begin{eqnarray}\label{e2.15}\nonumber
     \int_\Omega|u|^{p+1}&=&\frac{p+1}{p-1}\Bigg(2E(u)+\int_\Omega uu_t\Bigg)\\
       &\leq&\frac{2(p+1)}{p-1}E(u)+\frac{p+1}{p-1}\Bigg(\int_\Omega u^2\Bigg)^{1/2}\Bigg(\int_\Omega u_t^2\Bigg)^{1/2}.
   \end{eqnarray}
 Hereafter, we denote $C$ to be the constants varying from line to line. Using the monotonicity of energy and \eqref{e2.13}, after integrating \eqref{e2.15} over time, we obtain that
    \begin{eqnarray}\label{e2.16}\nonumber
      \int^\omega_0\int_\Omega|u|^{p+1}dxdt&\leq& C_p\Bigg(\omega E(u_0)+\omega+\int^\omega_0\int_\Omega u_t^2\Bigg)\\
       &\leq& C_p\Bigg(\omega E(u_0)+\omega+E(u_0)-{\mathcal{B}}\Bigg).
    \end{eqnarray}
 Now, integrating \eqref{e2.14} over $[0,\omega)$ and using \eqref{e2.16}, we get
    \begin{equation}\label{e2.17}
      \sup_{t\in[0,\omega)}\int_\Omega u^2(t)\leq 2\int^\omega_0\int_\Omega|u|^{p+1}+\int_\Omega u_0^2\leq C<\infty.
    \end{equation}
 Combining with \eqref{e2.15} and \eqref{e2.13}, it infers from \eqref{e2.17} that
   \begin{equation}\label{e2.18}
     \int^\omega_0\Bigg(\int_\Omega|u|^{p+1}dx\Bigg)^2dt\leq C.
   \end{equation}
 We claim that
   \begin{equation}\label{e2.19}
     \int^\omega_0\Bigg(\int_\Omega|u|^{p+1}dx\Bigg)^{\widetilde{q}}dt\leq C<\infty
   \end{equation}
 holds for all
   \begin{equation}\label{e2.20}
     \widetilde{q}<\max\Bigg\{q+\frac{2}{p+1}, \ \ \ \frac{2N}{N+2}\Bigg(q+\frac{2}{N+2}-\frac{2N}{(N+2)^2q+N(N+2)}\Bigg)\Bigg\},
   \end{equation}
 provided
   \begin{equation}\label{e2.21}
     \int^\omega_0\Bigg(\int_\Omega|u|^{p+1}dx\Bigg)^qdt\leq C<\infty
   \end{equation}
 for some $q\geq2$. The following argument is motivated by Quittner \cite{Q}, we present here for the conveniences of the readers. At first, it follows from \eqref{e2.13} \eqref{e2.21} and an interpolation theorem in \cite{CL} that
    \begin{equation}\label{e2.22}
      \sup_{t\in[0,\omega)}||u(t)||_{L^\lambda(\Omega)}\leq C<\infty
    \end{equation}
 for all
    $$
     \lambda<\lambda_1(q)\equiv p+1-\frac{p-1}{q+1}.
    $$
 Another hand, by Sobolev embedding theorem, we have
   \begin{eqnarray}\label{e2.23}\nonumber
     \int^\omega_0\Bigg(\int_\Omega|u|^{2^*}dx\Bigg)^{\frac{2q}{2^*}}dt&\leq& C\int^\omega_0\Bigg(\int_\Omega|\nabla u|^2dx\Bigg)^qdt\\
     &\leq& C\Bigg[\int^\omega_0\Bigg(\int_\Omega|u|^{p+1}dx\Bigg)^qdt+E(u_0)\Bigg]\leq C<\infty.
   \end{eqnarray}
 Using again the interpolation theorem in \cite{CL} together with \eqref{e2.13}\eqref{e2.23}, we conclude that \eqref{e2.22} holds for all
    $$
     \lambda<\lambda_2(q)\equiv\frac{2N(q+1)}{(N-2)q+N}.
    $$
To proceed further, applying the H\"{o}lder's inequality and Young's inequality to \eqref{e2.15}, we get
   \begin{eqnarray}\label{e2.24}\nonumber
     \int_\Omega|u|^{p+1}&\leq& C\Bigg\{E(u_0)+\Bigg(\int_\Omega |u|^\lambda\Bigg)^{\frac{1}{\lambda}}\Bigg(\int_\Omega|u_t|^{\frac{\lambda}{\lambda-1}}\Bigg)^{\frac{\lambda-1}{\lambda}}\Bigg\}\\ \nonumber
      &\leq& C\Bigg\{1+\Bigg(\int_\Omega|u_t|^{2\theta}|u_t|^{\frac{\lambda}{\lambda-1}-2\theta}\Bigg)^{\frac{\lambda-1}{\lambda}}\Bigg\}\\ \nonumber
      &\leq& C\Bigg\{1+\Bigg(\int_\Omega u_t^2\Bigg)^{\theta\frac{\lambda-1}{\lambda}}\Bigg(\int_\Omega|u_t|^{\frac{\lambda-2\theta(\lambda-1)}{(1-\theta)(\lambda-1)}}\Bigg)^{\frac{\lambda-1}{\lambda}(1-\theta)}\Bigg\}\\
      &\leq& C\Bigg\{1+\Bigg(\int_\Omega u_t^2\Bigg)^{\frac{1}{\widetilde{q}}}+\Bigg(\int_\Omega|u_t|^{\frac{\lambda-2\theta(\lambda-1)}{(1-\theta)(\lambda-1)}}\Bigg)^{\frac{(\lambda-1)(1-\theta)}{\lambda-\theta(\lambda-1)\widetilde{q}}}\Bigg\},
   \end{eqnarray}
 where $\theta, \widetilde{q}$ are chosen later. Noting that $|u|^{p-1}u\in L^{\frac{q(p+1)}{p}, \frac{p+1}{p}}(\Omega\times[0,\omega))$ and using $L^{\frac{q(p+1)}{p}, \frac{p+1}{p}}-$estimation for linear parabolic equation (see \cite{A}), we have $u_t\in  L^{\frac{q(p+1)}{p}, \frac{p+1}{p}}(\Omega\times[0,\omega))$. If $\lambda\geq p+1$, \eqref{e2.20} is clearly true since
    $$
     \sup_{t\in[0,\omega)}\int_\Omega(|\nabla u|^2+|u|^{p+1})(t)\leq C<\infty
    $$
 by \eqref{e2.22} and the boundedness of energy $E(u(t))$. So, we may assume that $2<\lambda<p+1$ due to $\lambda_2(q)>2$, and take
    $$
     \theta\equiv\frac{p+1-\lambda}{(p-1)(\lambda-1)}\in(0,1)
    $$
 such that
    \begin{equation}\label{e2.25}
     \frac{\lambda-2\theta(\lambda-1)}{(1-\theta)(\lambda-1)}=\frac{p+1}{p}.
    \end{equation}
 On the another hand, if we choose
   $$
    \widetilde{q}=\frac{q(p-1)\lambda}{(p-1)(\lambda-1)+(q-1)(p+1-\lambda)},
   $$
 then
   \begin{equation}\label{e2.26}
     \frac{(\lambda-1)(1-\theta)}{\lambda-\theta(\lambda-1)\widetilde{q}}\widetilde{q}=q.
   \end{equation}
 So, \eqref{e2.20} follows from \eqref{e2.25} and \eqref{e2.26} by integrating \eqref{e2.24} over time, and then sending $\lambda\uparrow\lambda_1(q)$, $\lambda\uparrow\lambda_2(q)$ separately. Now, a bootstrap argument yields
   \begin{equation}\label{e2.27}
     \int^\omega_0\Bigg(\int_\Omega|\nabla u|^2+\int_\Omega|u|^{p+1}\Bigg)^qdt\leq C_q<\infty
   \end{equation}
 for any $1\leq q<\infty$, where $C_q$ is a positive constant depending on $q$. Noting that
    $$
      \lim_{q\to\infty}\lambda_2(q)=\frac{2N}{N-2}>p+1
    $$
 for $p<p_S$, we obtain that
    \begin{equation}\label{e2.27}
       \sup_{t\in[0,\omega)}\int_\Omega|u|^{p+1}\leq C<\infty
    \end{equation}
 by \eqref{e2.22}. Since $p+1>\frac{N}{2}(p-1)$ when $p<p_S$, a well known result for semi-group \cite{W1} yields the boundedness of solution up to $\omega$. $\Box$\\

\begin{coro}\label{c2.1}
  Let $1<p<p_S$, all finite time blow-up are collapsing and hence complete.
\end{coro}

Corollary \ref{c2.1} is a direct consequence of Theorem \ref{t2.1} and Proposition \ref{p2.2}. It's also notable that a similar result has been proven in \cite{BC} for nonnegative solutions.

\vspace{40pt}

\section{Refined blowing-up II:  type II rates}

\begin{lemm}\label{l3.1} There exists a constant $C_{p,N}$ depending only on $p>1$ and $N\geq2$, such that
   \begin{equation}\label{e3.1}
     v(0,0)\leq C_{p,N}\Bigg(\fint_{P_R(0,0)}|v|^{p+1}dx\Bigg)^{\frac{1}{p+1}}
   \end{equation}
for any classical solution $v$ of
     \begin{equation}\label{e3.2}
        v_t-\triangle v\leq0\ \mbox{ in }\overline{P_R(0,0)}\equiv \overline{B_{\sqrt{\frac{N}{2\pi e}}R}}\times[-R^2,0],
     \end{equation}
where $\fint$ stands for the average integration.
\end{lemm}

\noindent\textbf{Proof.} By translating
   $$
     v(x,t)\rightarrow v(Rx,R^2t),
   $$
we need only prove the lemma for $R=1$. Noting first that the mean value property proved in \cite{E} for heat equation still holds for \eqref{e3.2}, no other than replacing equalities by inequalities, we have
  \begin{equation}\label{e3.3}
    v(0,0)\leq\frac{1}{4r^N}\iint_{\widetilde{P}_r(0,0)}v(y,s)\frac{|y|^2}{s^2}dyds, \ \forall 0<r\leq R=1,
  \end{equation}
for
  $$
    \widetilde{P}_r(0,0)\equiv\Big\{(y,s)\in{\mathbb{R}}^{N+1}:\ s\leq0, \Phi(-y,-s)\geq\frac{1}{r^N}\Big\}\subset B_{\sqrt{\frac{N}{2\pi e}}r}\times\Big(-\frac{r^2}{4\pi},0\Big)\subset P_r(0,0),
  $$
where
  $$
    \Phi(x,t)=\frac{1}{(4\pi t)^{\frac{N}{2}}}e^{-\frac{|x|^2}{4t}}, \ t>0, x\in{\mathbb{R}}^N
  $$
is the standard backward heat kernel. Since for any $1<q<\frac{N+2}{2}$, there exists a positive constant $C_{q,N}$ depending on $q$ and $N$, such that
     \begin{equation}\label{e3.4}
       \Bigg(\iint_{\widetilde{P}_1(0,0)}\frac{|y|^{2q}}{|s|^{2q}}\Bigg)^{\frac{1}{q}}\leq C_{q,N},
     \end{equation}
  the conclusion follows from \eqref{e3.3} and \eqref{e3.4} together with H\"{o}lder's inequality. $\Box$\\

\begin{lemm}\label{l3.2}
  Let $p>1$ and $u$ be a maximal classical solution to \eqref{e1.1} with maximal life time $\omega<\infty$. There exists a positive constant $\varepsilon_0$ depending only $p$ and $N$, such that if
    $$
      r^{\frac{4}{p-1}-N}\iint_{P_r(\overline{z})}|u|^{p+1}dxdt<\varepsilon_0
    $$
  holds for all cylinders $P_r(\overline{z})\equiv B_{\sqrt{\frac{N}{2\pi e}}r}(\overline{x})\times(\overline{t}-r^2,\overline{t})$ contained inside the cylinder $P_{R}(a,\omega)$, then $(a,\omega)$ is not a singular point. Moreover, there holds
     $$
       \sup_{P_\frac{R}{2}(a,\omega)}|u|\leq 4^{\frac{2}{p-1}}R^{-\frac{2}{p-1}}.
     $$
\end{lemm}

\noindent\textbf{Proof.} Without loss of generality, we may assume that $u(x,t)$ is smooth up to singular time $t=\omega$. Otherwise, we can shift $\omega$ to $\omega-\frac{1}{i}$ and then send $i\to+\infty$.

 Consider
   $$
     M=\sup_{0<r<R}\Big[(R-r)^{\frac{2}{p-1}}\sup_{\overline{P_r(a,\omega)}}|u|\Big]
   $$
 and let $r_0\in[0,R), z^*\in\overline{P_{r_0}(a,\omega)}$ such that
   $$
     M=(R-r_0)^{\frac{2}{p-1}}|u|(z^*).
   $$
 Setting $r_1=\frac{R-r_0}{2}$, we have
   $$
     P_{r_1}(z^*)\subset P_R(a,\omega)
   $$
 and
    $$
      r_1^{\frac{2}{p-1}}\sup_{P_{r_1}(z^*)}|u|\leq M.
    $$
 It yields that
   $$
     \sup_{P_{r_1}(z^*)}|u|\leq\Big(\frac{R-r_0}{r_1}\Big)^{\frac{2}{p-1}}|u|(z^*)=4^{\frac{1}{p-1}}|u|(z^*).
   $$
 Re-scaling $u$ by
   $$
     \widetilde{u}(y,s)=\mu^{-1}u\Big(x^*+\mu^{\frac{1-p}{2}}y,t^*+\mu^{1-p}s\Big),\ \ \ \mu\equiv u(z^*),
   $$
 then $\widetilde{u}$ satisfies that
   \begin{equation}\label{e3.5}
     \begin{cases}
       \widetilde{u}_s=\triangle\widetilde{u}+|\widetilde{u}|^{p-1}\widetilde{u},\\
       |\widetilde{u}|\leq4^{\frac{1}{p-1}},\ |\widetilde{u}|(0,0)=1
     \end{cases}
   \end{equation}
 in $P_{\mu^{\frac{p-1}{2}}r_1}(0,0)$. We claim that $M\leq 4^{\frac{1}{p-1}}$. For, if $M>4^{\frac{1}{p-1}}$, then
   $$
     \mu^{\frac{p-1}{2}}r_1\geq1
   $$
 and \eqref{e3.5} holds in $P_1(0,0)$. Noting that
    \begin{equation}\label{3.6}
      \iint_{P_1(0,0)}|\widetilde{u}|^{p+1}dyds=\mu^{\frac{p-1}{2}N-2}\iint_{P_{\mu^{\frac{1-p}{2}}}(z^*)}|u|^{p+1}dxdt<\varepsilon_0
    \end{equation}
 and regarding \eqref{e3.5} as a linear parabolic inequality
   $$
     \Big(\frac{\partial}{\partial s}-\triangle\Big)\widetilde{u}^2=-2|\nabla\widetilde{u}|^2+2|\widetilde{u}|^{p+1}\leq 8\widetilde{u}^2,
   $$
 after applying Lemma \ref{l3.1} to
   $$
     v(y,s)=e^{-8s}\widetilde{u}^2(y,s),
   $$
 we get
   \begin{eqnarray*}
     1=\widetilde{u}^2(0,0)&=&|v|(0,0)\\
       &\leq& C_{p,N}\Big(\fint_{P_1(0,0)}|v|^{p+1}dyds\Big)^{\frac{1}{p+1}}\\
       &\leq& C'_{p,N}\Big(\fint_{P_1(0,0)}|\widetilde{u}|^{2p+2}dyds\Big)^{\frac{1}{p+1}}\\
       &\leq& C''_{p,N}\varepsilon_0^{\frac{1}{p+1}}.
   \end{eqnarray*}
 Contradiction holds provided $\varepsilon_0$ is chosen small. So $M\leq 4^{\frac{1}{p-1}}$ and hence
   $$
     \Big(\frac{R}{2}\Big)^{\frac{2}{p-1}}\sup_{P_{\frac{R}{2}}(a,\omega)}|u|\leq 4^{\frac{1}{p-1}}.
   $$
 The conclusion is drawn. $\Box$\\

Under below, we transform the maximal solution to self-similar variables
   $$
     u(x,t)=(\omega-t)^{-\frac{1}{p-1}}w\Bigg(\frac{x-a}{\sqrt{\omega-t}},-\log(\omega-t)\Bigg)
   $$
which was inspired by Giga-Kohn \cite{GK1}. The re-scaled function $w=w_{(a,\omega)}(y,s)$ satisfies that
  \begin{equation}\label{e3.7}
    w_s-\triangle w+\frac{1}{2}y\cdot\nabla w+\frac{1}{p-1}w=|w|^{p-1}w
  \end{equation}
in $\cup_{s>-\log\omega}\Omega_s$ for
   $$
     \Omega_s=\Omega_s^a\equiv\{y\in{\mathbb{R}}^N:\ e^{-\frac{s}{2}}y+a\in\Omega\}.
   $$
Re-writing \eqref{e3.7} as self-adjoint form
  \begin{equation}\label{e3.8}
    \rho w_s-\frac{\partial}{\partial y_i}\Big(\rho\frac{\partial w}{\partial y_i}\Big)+\frac{1}{p-1}\rho w=\rho|w|^{p-1}w
  \end{equation}
with $\rho(y)=e^{-\frac{|y|^2}{4}}$, we derive that
 \begin{equation}\label{e3.9}
   \frac{d}{ds}{\mathcal{E}}(w(s))=-\int_{\Omega_s}w_s^2\rho dy-\frac{1}{4}\int_{\partial\Omega_s}(y\cdot\nu)\Big|\frac{\partial w}{\partial\nu}\Big|^2\rho d\sigma
 \end{equation}
and
  \begin{equation}\label{e3.10}
    \frac{1}{2}\frac{d}{ds}\int_{\Omega_s}w^2\rho dy=-2{\mathcal{E}}(w)+\frac{p-1}{p+1}\int_{\Omega_s}|w|^{p+1}\rho dy,
  \end{equation}
after multiplying by $w_s$, $w$ respectively and then integrating over $\Omega_s$, where
   $$
     {\mathcal{E}}(w(s))\equiv\frac{1}{2}\int_{\Omega_s}|\nabla w|^2\rho dy+\frac{1}{2(p-1)}\int_{\Omega_s}w^2\rho dy-\frac{1}{p+1}\int_{\Omega_s}|w|^{p+1}\rho dy
   $$
is the local energy of $w$.

\vspace{20pt}

The following characteristic property of blowup expressed in term of local energy ${\mathcal{E}}(w)$, was firstly proved in \cite{GK2} under subcritical case:

\begin{theo}\label{t3.1}
  Let $\Omega$ be convex and $u$ be a maximal classical solution to \eqref{e1.1} with $p>1$. There exists a positive constant $\varepsilon_1$ depending only on $p$ and $N$, such that
    $$
     {\mathcal{E}}(w_{a,\omega}(s_0))\leq\varepsilon_1\ \ \mbox{ for some } s_0>-\log\omega
    $$
  implies $(a,\omega)$ is not a singular point.
\end{theo}

\noindent\textbf{Remark 3.1} It's not hard to remove the convexity assumption of the theorem, by addapting a similar argument as in \cite{CD}.\\

\noindent\textbf{Proof.} We will prove that for some $\delta_0>0$ (depending on $u,s_0$), there holds
  \begin{equation}\label{e3.11}
    ess\sup_{P_{\frac{\delta_0}{2}}(a,\omega)}|u|\leq 4^{\frac{2}{p-1}}\delta_0^{-\frac{2}{p-1}}.
  \end{equation}

 By continuity, $\exists 0<\delta_0<\frac{1}{2}e^{-\frac{s_0}{2}}$, such that for any
   $$
     (a',\omega')\in\overline{B_{\sqrt{\frac{N}{2\pi e}}\delta_0}(a)\times(\omega-\delta_0^2,\omega+\delta_0^2)},
   $$
 there holds
  \begin{equation}\label{e3.12}
     {\mathcal{E}}(w_{(a',\omega')}(s'_0))\leq 2\varepsilon_1,
  \end{equation}
where
 $$
    s'_0\equiv-\log(\omega'-\omega+e^{-s_0})\in\Big[s_0-\log\frac{5}{4},s_0+\log\frac{4}{3}\Big].
 $$
Applying \eqref{e3.9} to $w=w_{(a',\omega')}(y,s)$, a function defined over $[s'_0, -\log(\omega'-\omega)_+$), we have ${\mathcal{E}}(w(s))$ is a monotone non-increasing function on $s$. Thus,
  \begin{equation}\label{e3.13}
    \int^{s+1}_s\int_{\Omega_\tau}w_s^2\rho dyd\tau\leq{\mathcal{E}}(w(s))\leq{\mathcal{E}}(w(s'_0))
  \end{equation}
for all $s'_0\leq s<-\log(\omega'-\omega)_+-1$. Now, taking a small constant $\varepsilon'_1>0$ and using Young's inequality, we conclude that
  \begin{eqnarray*}
    \int^{s+1}_s\int_{\Omega_\tau}|w|^{p+1}\rho dyd\tau&\leq&\frac{2(p+1)}{p-1}\int^{s+1}_s{\mathcal{E}}(w(\tau))d\tau+\frac{p+1}{p-1}\int^{s+1}_s\int_{\Omega_\tau}ww_s\tau\rho dyd\tau\\
    &\leq& \frac{2(p+1)}{p-1}{\mathcal{E}}(w(s'_0))+\frac{1}{2}\int^{s+1}_s\int_{\Omega_\tau}|w|^{p+1}\rho dyd\tau+\frac{\varepsilon'_1}{4}\\
    && + C_p{\varepsilon'_1}^{-\frac{p-1}{p+1}}\int^{s+1}_s\int_{\Omega_\tau}w_s^2\rho dyd\tau\\
    &\leq&\Bigg(\frac{2(p+1)}{p-1}+C_p{\varepsilon'_1}^{-\frac{p-1}{p+1}}\Bigg){\mathcal{E}}(w(s'_0))+\frac{1}{2}\int^{s+1}_s\int_{\Omega_\tau}|w|^{p+1}\rho dyd\tau+\frac{\varepsilon'_1}{4}
  \end{eqnarray*}
by \eqref{e3.10}, where
   $$
     C_p=\frac{p+1}{2}\Big(2(p-1)\Big)^{-\frac{p+3}{p+1}}.
   $$
Consequently, we obtain that for all $s'_0\leq s<-\log(\omega'-\omega)_+-1$,
  \begin{equation}\label{e3.14}
    \int^{s+1}_s\int_{\Omega_\tau}|w|^{p+1}\rho dyd\tau\leq\Bigg(\frac{4(p+1)}{p-1}+2C_p{\varepsilon'_1}^{-\frac{p-1}{p+1}}\Bigg){\mathcal{E}}(w(s'_0))+\frac{\varepsilon'_1}{2}\leq\varepsilon'_1,
  \end{equation}
as long as $\varepsilon_1$ is small.

 Writing \eqref{e3.14} back to $u(x,t)$, we get
   \begin{equation}\label{e3.15}
      r^{\frac{4}{p-1}-N}\int^{\omega'-\frac{r^2}{e}}_{\omega'-r^2}\int_{B_{\sqrt{\frac{N}{2\pi e}}r}(a')}|u|^{p+1}dxdt\leq\varepsilon''_1\equiv e^{\frac{2}{p-1}+\frac{N}{8\pi e}-\frac{N}{2}}\varepsilon'_1
   \end{equation}
 for any $0<r<\frac{4}{3}e^{-\frac{s_0}{2}}$. Therefore, by changing the parameters
   $$
    \overline{x}=a',\ \overline{t}=\omega'-\frac{r^2}{e},\ \overline{r}=\sqrt{\frac{e-1}{e}}r
   $$
 it's inferred from \eqref{e3.15} that
  \begin{equation}\label{e3.16}
    r^{\frac{4}{p-1}-N}\iint_{P_{\overline{r}}(\overline{z})}|u|^{p+1}dxdt\leq\varepsilon''_1<\varepsilon_0
  \end{equation}
 for any $P_{\overline{r}}(\overline{z})$ contained inside $P_{\delta_0}(a,\omega)$, suppose that $\varepsilon_1$ is chosen small. Now, the conclusion follows from Lemma \ref{l3.2}. $\Box$\\

 Next, let's recall a key decaying estimation which was proven by Giga-Kohn (\cite{GK1}, Proposition 4.3 and Remark 4.4) under critical and subcritical case:

 \begin{prop}\label{p3.1}
    Assume that $\Omega$ is star-shaped with respect to $a\in\Omega$, and let $w=w_{(a,\omega)}(y,s)$ be the re-scaled solution of \eqref{e3.7}. There exists a positive constant $C_*$ which is independent of $a$, such that
      $$
        \int^\infty_{-\log\omega}\int_{\Omega_s}\rho(1+|y|^2)(w_s^2+|\nabla w|^2)dyds\leq C_*<\infty
      $$
 holds in case of $1<p<p_S$, and
      $$
        \int^\infty_{-\log\omega}\int_{\Omega_s}\rho(1+|y|^2)w_s^2+\rho|y|^2|\nabla w|^2dyds\leq C_*<\infty.
      $$
 holds in case of $p=p_S$.
 \end{prop}

To control the quadratic term in local energy, we can argue as Giga-Kohn to get a slightly stronger version of Lemma 4.1 in \cite{GK2}:

\begin{lemm}\label{l3.3}
  There exists a positive constant $C_N$ depending only on $N$, such that for any $H^1_{loc}$ function $f$ on ${\mathbb{R}}^N$ and $\alpha>0$,
    $$
      \int_{|y|>\alpha}\rho f^2dy\leq C_N\Bigg\{\int_{|y|>\alpha}\rho\Big(|y|^2+|y|^{-(N-1)}\Big)|\nabla f|^2dy+\Big(\fint_{|y|=\alpha}fd\sigma\Big)^2\Bigg\},
    $$
  where $\rho(y)=e^{-\frac{|y|^2}{4}}$ and $d\sigma$ stands for the area element of the sphere.\\
\end{lemm}

\noindent Now, taking a positive non-increasing function $\alpha(s)$ for $s\in{\mathbb{R}}$, we assume that

(H1) In case of $\lim_{s\to+\infty}\alpha(s)=0$, there exists some positive constant $\alpha_0$ such that $\alpha(s)\geq \alpha_0a^{-\frac{1}{N+1}}(s)$ for all $s\geq1$, where $a(\cdot)$ is a positive non-decreasing function satisfying
     \begin{equation}\label{e3.17}
        \int^{+\infty}_0\frac{ds}{a(s)}=+\infty.
     \end{equation}

(H2)  In case of $\lim_{s\to+\infty}\alpha(s)>0$, there holds $\alpha(s)\geq\widetilde{\alpha}_0$ for some $\widetilde{\alpha}_0$ and all $s\geq1$.\\

Next, we can estimate the local energy by separating the integral region into inner and outer ones as following:

\begin{prop}\label{p3.2}
   Assume that $\Omega$ is star-shaped with respect to $a\in\Omega$ and let $w=w_{(a,\omega)}(y,s)$ be the re-scaled solution of \eqref{e3.7}. Taking a function $\alpha(s)$ fulfilling (H1) or (H2), we have the following estimations on local energy of $w=w_{(a,\omega)}(y,s)$:

 \noindent Case 1: $\lim_{s\to+\infty}\alpha(s)=0$, there holds
   \begin{eqnarray}\label{e3.18}\nonumber
      \int^{s+1}_s{\mathcal{E}}(w(\tau))d\tau&\leq&\alpha_0^{-N-1}C_N\Bigg(a(s)\int^{s+1}_s\int_{\Omega_\tau}|y|^2\rho|\nabla w|^2dyd\tau\Bigg)+C_{p,N}\int^{s+1}_s\alpha^N(\tau)d\tau\\
       &&+\frac{1}{2}\int^{s+1}_s\int_{|y|\leq\alpha(\tau)}|\nabla w|^2dyd\tau+C_N\int^{s+1}_s\Big(\fint_{|y|=\alpha(\tau)}wd\sigma\Big)^2d\tau,
   \end{eqnarray}

\noindent Case 2: $\liminf_{s\to+\infty}\alpha(s)>0$, there holds
   \begin{eqnarray}\label{e3.19}\nonumber
      \int^{s+1}_s{\mathcal{E}}(w(\tau))d\tau&\leq&\widetilde{\alpha}_0^{-N-1}C_N\int^{s+1}_s\int_{\Omega_\tau}|y|^2\rho|\nabla w|^2dyd\tau+C_N\int^{s+1}_s\int_{|y|\leq\alpha(\tau)}|\nabla w|^2dyd\tau\\
       &&+C_N\Bigg\{\int^{s+1}_s\Big(\fint_{|y|=\alpha(\tau)}wd\sigma\Big)^2d\tau
       +\int^{s+1}_s\Big(\fint_{|y|\leq\alpha(\tau)}wdy\Big)^2d\tau\Bigg\},
   \end{eqnarray}
where $C_N>0$ is a constant depending only on $N$ and $C_{p,N}>0$ is a constant depending only on $p,N$.
\end{prop}

It's notable that the first two terms on right hand side of \eqref{e3.18} and the first term on right hand side of \eqref{e3.19} are all tends to zero for a subsequence $s=s_k\to\infty$ by Proposition \ref{p3.1}. Therefore, controlling of the remain terms will force the local energy of $w$ to dissipate, and so exclude the possibility of blowing-up.\\

\noindent\textbf{Proof of Proposition \ref{p3.2}:} The proof relies on the separation of the integral region from $\{|y|>\alpha\}$ and $\{|y|\leq\alpha\}$ for $\alpha=\alpha(s)$. In fact, we have
  \begin{equation}\label{e3.20}
    \int^{s+1}_s\int_{|y|>\alpha(\tau)}\rho(y)|\nabla w|^2dyd\tau\leq\frac{a(s)}{a(s)\alpha^2(s+1)}\int^{s+1}_s\int_{|y|>\alpha(\tau)}|y|^2\rho(y)|\nabla w|^2dyd\tau,
  \end{equation}
and
  \begin{eqnarray}\label{e3.21}\nonumber
    \int^{s+1}_s\int_{|y|>\alpha(\tau)}\rho w^2dyd\tau&\leq&C_N\frac{a(s)}{a(s)\alpha^{N+1}(s+1)}\int^{s+1}_s\int_{|y|>\alpha(\tau)}|y|^2\rho|\nabla w|^2dyd\tau\\
    &&+C_N\int^{s+1}_s\Bigg(\fint_{|y|=\alpha(\tau)}wd\sigma\Bigg)^2d\tau
  \end{eqnarray}
after using Lemma \ref{l3.3} for $f=w$ and integrating over time. Therefore, the inner region can be controlled as following: If $\lim_{s\to+\infty}\alpha(s)=0$, we estimate
  \begin{eqnarray}\label{e3.22}\nonumber
    \frac{1}{2(p-1)}\int^{s+1}_s\int_{|y|\leq\alpha(\tau)}\rho w^2dyd\tau&\leq&\frac{p-1}{2(p+1)}\int^{s+1}_s\int_{|y|\leq\alpha(\tau)}\rho|w|^{p+1}dyd\tau\\
    &&+C_{p,N}\int^{s+1}_s\alpha^N(\tau)d\tau
  \end{eqnarray}
by Young's inequality. If $\liminf_{s\to+\infty}\alpha(s)>0$, we estimate
  \begin{equation}\label{e3.23}
    \frac{1}{2(p-1)}\int^{s+1}_s\int_{|y|\leq\alpha(\tau)}\rho w^2dyd\tau\leq C_N\int^{s+1}_s\int_{|y|\leq\alpha(\tau)}|\nabla w|^2+C_N\int^{s+1}_s\Bigg(\fint_{|y|\leq\alpha(\tau)}wdy\Bigg)^2d\tau
  \end{equation}
after using Poincar\'{e}'s inequality
  \begin{eqnarray}\label{e3.24}\nonumber
    \int_{|y|\leq\alpha(\tau)}w^2dy&\leq&2\int_{|y|\leq\alpha(\tau)}(w-\overline{w})^2dy+2C_N\alpha^N(\tau)\overline{w}^2\\ &\leq&C_N\alpha^2(\tau)\int_{|y|\leq\alpha(\tau)}|\nabla w|^2+C_N\alpha^N(\tau)\Bigg(\fint_{|y|\leq\alpha(\tau)}w\Bigg)^2
  \end{eqnarray}
and then integrating over time, where $\overline{w}$ stands for the average of $w$ on $B_{\alpha(\tau)}$. So, a combination of \eqref{e3.20} to \eqref{e3.23} yields the desired conclusion. $\Box$\\

Combining the decaying estimation Proposition \ref{p3.1} with Proposition \ref{p3.2}, we have the following property:

\begin{prop}\label{p3.3}
   Assume that $\Omega$ is convex and let $w=w_{(a,\omega)}(y,s)$ be the re-scaled solution of maximal solution $u(x,t)$ to \eqref{e1.1} over $[0,\omega)$, where $1<p\leq p_s$. There exist a positive constant $\varepsilon_2$ depending only on $p$ and $N$, such that $(a,\omega)$ is not a singular point for $u$, provided

Case1: when $\alpha(s)$ fulfills (H1) with $\lim_{s\to+\infty}\alpha(s)=0$, assume further that
\begin{equation}\label{e3.25}
    \limsup_{s\to+\infty}\Bigg\{\int^{s+1}_s\int_{|y|\leq\alpha(\tau)}|\nabla w|^2dyd\tau+\int^{s+1}_s\Big(\fint_{|y|=\alpha(\tau)}wd\sigma\Big)^2d\tau\Bigg\}<\varepsilon_2,
  \end{equation}
or

Case2: when $\alpha(s)$ fulfills (H2) with $\liminf_{s\to+\infty}\alpha(s)>0$, assume further that
\begin{eqnarray}\label{e3.26}\nonumber
    &&\inf_{s\geq s_0}\Bigg\{\int^{s+1}_s\int_{|y|\leq\alpha(\tau)}|\nabla w|^2dyd\tau+\int^{s+1}_s\Big(\fint_{|y|=\alpha(\tau)}wd\sigma\Big)^2d\tau\\
    && \ \ \ \ \ \ \ \ \ \ \ \ \  \ \ \ \ \ \ \ \ \ \ \ \ \  \ \ \ \ \ \ \ \ \ \ \ \ \  \ \ \ \ \ \ \ \ \ \ \ \ \ +\int^{s+1}_s\Big(\fint_{|y|\leq\alpha(\tau)}wdy\Big)^2d\tau\Bigg\}<\varepsilon_2
  \end{eqnarray}
for some large constant $s_0$ depending only $p, N$ and $u_0$.
\end{prop}

\noindent\textbf{Proof.} The key observation is that
  $$
   \liminf_{s\to+\infty}a(s)\int^{s+1}_s\int_{\Omega_\tau}|y|^2\rho|\nabla w|^2dyd\tau=0
  $$
or
  $$
   \lim_{s\to+\infty}\int^{s+1}_s\int_{\Omega_\tau}|y|^2\rho|\nabla w|^2dyd\tau=0
  $$
owing to Proposition \ref{p3.1}. Combined with Proposition \ref{p3.2}, we get
  $$
    \liminf_{s\to+\infty}\int^{s+1}_s{\mathcal{E}}(w(\tau))d\tau<\varepsilon_1,
  $$
and hence
   $$
    \lim_{s\to+\infty}{\mathcal{E}}(w(s))<\varepsilon_1
   $$
by intermediate value theorem for integral and monotonicity of the local energy. So, conclusion follows from Theorem \ref{t3.1}. $\Box$\\

\noindent\textbf{Theorem 3.2} (Non-degeneracy of blowing-up)
  Under assumptions of Proposition 3.3 and assume that $p=p_S$, we take a function $\alpha(s)\equiv1$. Then $(a,\omega)$ is not a blowing-up point for $u$, provided
    \begin{equation}\label{e0.1}
      \limsup_{t\to\omega^-}\big(2\sqrt{\omega-t}\big)^{\frac{2}{p-1}}\sup_{\{x\in\Omega:\ |x-a|\leq 2\sqrt{\omega-t}\}}|u|(x,t)=0.
    \end{equation}

\noindent\textbf{Proof.} Under the self-similar variables, it yields from \eqref{e0.1} that
  \begin{equation}\label{e0.2}
    \limsup_{s\to+\infty}2^{\frac{2}{p-1}}\sup_{\{y\in\Omega_s:\ |y|\leq2\}}|w|(y,s)=0.
  \end{equation}
Noting that $w$ satisfies that
  \begin{equation}\label{e0.3}
    w_s-\triangle w+\frac{1}{2}y\cdot\nabla w+\frac{1}{p-1}w=|w|^{p-1}w
  \end{equation}
in $B_2\times\big[-\frac{2}{\lambda^2},0\big]$ and $w$ becomes small as $s$ being large. Regarding \eqref{e0.3} as a linear parabolic equation with bounded coefficients, we conclude that $\sup_{B_1\times[-\frac{1}{\lambda^2},0]}|\nabla w|$ also tends to zero as $s$ tends to infinity. We conclude that
   \begin{eqnarray}\label{e0.4}\nonumber
     &&\int^{s+1}_s\int_{|y|\leq 1}|\nabla w|^2dyd\tau+\int^{s+1}_s\Big(\fint_{|y|=1}wd\sigma\Big)^2d\tau\\
    && \ \ \ \ \ \ \ \ \ \ \ \ \  \ \ \ \ \ \ \ \ \ \ \ \ \  \ \ \ +\int^{s+1}_s\Big(\fint_{|y|\leq 1}wdy\Big)^2d\tau\to0
   \end{eqnarray}
 when $s\to+\infty$. As a result, the conclusion follows from Proposition (3.3). $\Box$\\

Next, we will focus on non-collapsing blowup, and prove that it must be type II in sense of (1.3). It notable that for non-collapsing blowup, we have
   \begin{equation}\label{e0.5}
     \int^\omega_0\int_\Omega u_t^2dxdt<+\infty.
   \end{equation}

 Now, fixing $a\in\Omega, b\in(0,1]$, we denote a shrinking rotary paraboloid by
    $$
      Q_{(a,\omega)}\equiv\Big\{z=(a+y\sqrt{\omega-t},t)\in\Omega\times[0,\omega):\ |y|\leq 1, 0\leq t\leq \omega\Big\}
    $$
of width $b$ and centering at $(a,\omega)$, and denote its early part by
    $$
      Q_{(a,\omega)}(t_0)\equiv\Big\{z=(a+y\sqrt{\omega-t},t)\in\Omega\times[0,\omega):\ |y|\leq 1, 0\leq t\leq t_0\Big\}
    $$
for very $0\leq t_0\leq \omega$.

  Given any positive number $M>\max_{\Omega}|u_0|$, we set
     \begin{equation}\label{e0.6}
       \Lambda(M)\equiv\Lambda_{(a,\omega)}(M)=\inf\Big\{0\leq t'<\omega:\ \sup_{\overline{z}=(\overline{x},\overline{t})\in Q_{(a,\omega)}(t')}|u|(\overline{z})=M\Big\}
     \end{equation}
  and define
     \begin{equation}\label{e0.7}
       \lambda(t)\equiv\lambda_{(a,\omega)}(t)=\sup_{\overline{z}=(\overline{x},\overline{t})\in Q_{(a,\omega)}(t)}|u|(\overline{z})
     \end{equation}
  for any $0<t<\omega$.\\

we can now state the following result of type II blowup for energy non-collapsing in critical case:\\

\noindent\textbf{Theorem 3.3}  Assume $\Omega$ is convex and $p=p_s, n\geq3$. Let $u$ be a maximal solution which blows up at $\omega<+\infty$. Assume that the energy is non-collapsing, then all blowing-up points $(a,\omega)$ must be of type II.  More precisely, we have
      \begin{equation}\label{e0.8}
        \limsup_{t\to \omega^-}(\omega-t)^{\frac{1}{p-1}}\lambda_{(a,\omega)}(t)=+\infty.
      \end{equation}

\vspace{5pt}

It's notable that the trivial type I singular solution
   $$
     u(x,t)=(p-1)^{-\frac{1}{p-1}}(\omega-t)^{-\frac{1}{p-1}}
   $$
has been ruled out since it blows up collapsing. Let's start with a corollary of Theorem 3.2:\\

\noindent\textbf{Lemma 3.4}  Assume that $\Omega$ is convex and $p=p_s, n\geq3$. Let $u(x,t)$ be a maximal solution which blows up at $\omega<+\infty$. Then for any blowing-up point $a\in\overline{\Omega}$, we have
       \begin{equation}\label{e0.9}
         \liminf_{t\to\omega^-}(\omega-t)^{\frac{1}{p-1}}\lambda_{(a,\omega)}(t)>0
       \end{equation}

\noindent\textbf{Remark 3.2} it follows from Lemma 3.4 that
   \begin{equation}\label{e0.10}
     \lim_{t\to\omega^-}\lambda_{(a,\omega)}(t)=+\infty.
   \end{equation}

\vspace{10pt}

\noindent\textbf{Proof of Theorem 3.3} We will prove \eqref{e0.8} by contradiction. Suppose on the contrary, there must be a large constant $C^*>0$ such that
  \begin{equation}\label{e0.11}
    \frac{1}{C^*}\leq(\omega-t)^{\frac{1}{p-1}}\lambda_{(a,\omega)}(t)\leq C^*
  \end{equation}
for all $0<t<\omega$ due to Lemma 3.4. We claim now there exists a positive constant $C_*$ such that
  \begin{equation}\label{e0.12}
    \liminf_{M\to+\infty}(M)^{p-1}(\Lambda(2M)-\Lambda(M))\leq C_*<+\infty.
  \end{equation}
For, if not, then
   \begin{equation}\label{e0.13}
    \lim_{M\to+\infty}M^{p-1}(\Lambda(2M)-\Lambda(M))=+\infty.
   \end{equation}
Consequently,
   $$
    \varphi(k)\equiv(2^{k})^{p-1}[\Lambda(2^{k+1})-\Lambda(2^{k})]\to+\infty   $$
as $k\to+\infty$. Assume $0<t<\omega$ such that
 $$
   2^{k}<\lambda(t)\leq 2^{k+1}
 $$
for some $k\in{\mathbb{N}}$, we have
  $$
   \Lambda(2^{k})<t\leq \Lambda(2^{k+1})
  $$
by definition. Thus, one can prove
   $$
     (\omega-t)^{\frac{1}{p-1}}\lambda(t)\geq\Big(\Lambda(2^{k+2})-\Lambda(2^{k+1})\Big)^{\frac{1}{p-1}}\cdot 2^{k}=\frac{1}{2}\varphi^{\frac{1}{p-1}}(k+1)\to+\infty,
   $$
 which contradicts with \eqref{e0.11}. Therefore, \eqref{e0.12} holds true.

Now, take a monotone increasing sequence $\{M_k\}_{k=1}^\infty$ which tends to infinity as $k\to\infty$, and such that
     \begin{equation}\label{e0.14}
      M_k^{p-1}[\Lambda(2M_k)-\Lambda(M_k)]\leq C_*<+\infty
     \end{equation}
  for all $k\in{\mathbb{N}}$. Furthermore,
     \begin{equation}\label{e0.15}
       \Lambda(M_k)<\omega
     \end{equation}
  for all $k\in{\mathbb{N}}$ by \eqref{e0.10}.

   By definition, there exist some sequences $x_k\in\Omega$
   and $t_k=\Lambda(2M_k)$, such that
      \begin{equation}\label{e0.16}
        |u(x_k,t_k)|=\sup_{(x,t)\in Q_{(a,\omega)}(t_k)}|u|(x,t)=2M_k.
      \end{equation}
  Now, re-scaling $u$ by
     $$
      u_k(x,t)=\frac{1}{2M_k}u\Big(\frac{x}{(2M_k)^{\frac{p-1}{2}}}+x_k,
      \frac{t}{(2M_k)^{p-1}}+t_k\Big),
     $$
  we have the re-scaled function defined on
    \begin{eqnarray*}
     Q^k=\Big\{(x,t)\in\Omega\times[0,\omega):\ |x-a_k|\leq b_k(t), t\in[-t_k(2M_k)^{p-1},0]\Big\},
    \end{eqnarray*}
where
  $$
    a_k\equiv(a-x_k)\cdot(2M_k)^{\frac{p-1}{2}},
  $$
  $$
    b_k(t)\equiv \sqrt{\omega-t}\cdot(2M_k)^{\frac{p-1}{2}}.
  $$
Noting that
  $$
   |a_k|\leq b_k(t_k)=\sqrt{\omega-t_k}(2M_k)^{\frac{p-1}{2}}\leq (C^*)^{\frac{p-1}{2}}
  $$
by \eqref{e0.11}, for some subsequence $k=k_j$, we get limits
  $$
   a_\infty=\lim_{j\to+\infty}a_{k_j}, \ b_\infty=\lim_{j\to+\infty}b_{k_j}(t_{k_j}).
  $$
Another hand, it follows from \eqref{e0.11} that
  $$
   b_k(t)\geq b_k(t_k)=\sqrt{\omega-t_k}(2M_k)^{\frac{p-1}{2}}\geq (C^*)^{-\frac{p-1}{2}}.
  $$
So, $b_\infty\geq (C^*)^{-\frac{p-1}{2}}$. Consequently,
  \begin{equation}\label{e0.17}
    Q^k \supseteq B_{b_k(t_k)}(a_k)\times\Big[-t_k(2M_k)^{p-1},0\Big]\to B_{b_\infty}(a_\infty)\times(-\infty,0]\ni(0,0)
  \end{equation}
for $k=k_j$ and $j$ large. Furthermore,
   $$
     |u_k|(0,0)=\sup_{Q^k}|u_k|(x,t)=1
   $$
and
   $$
    |u_k|(x,t)\leq\frac{1}{2}
   $$
for all $(x,t)\in B_{b_k(t_k)}(a_k)\times[-t_k(2M_k)^{p-1},-N\cdot 2^{p-1}]$.

  Passing $k\to\infty$, we get a ancient solution $v(x,t)$
  defining on $Q^\infty\equiv B_{b_\infty}(a_\infty)\times(-\infty,0)\ni(0,0)$, and satisfying
    \begin{equation}\label{e0.18}
      \begin{cases}
        v_t-\triangle v=|v|^{p-1}v\ \ \  \mbox{ in } Q^\infty\\[5pt]
       |v|(0,0)=\sup_{Q^\infty}|v|(x,t)=1
      \end{cases}
    \end{equation}
  and
     \begin{equation}\label{e0.19}
       |v|(x,t)\leq\frac{1}{2}\ \ \ \forall (x,t)\in B_{b_\infty}(a_\infty)\times(-\infty,-N\cdot
   2^{p-1}).
     \end{equation}

 Another hand, we have
 \eqref{e0.5} for non-collapsing blowup.  Thus, for any fixed $R>0$, we have
    \begin{eqnarray*}
     \iint_{B_{b_k}(a_k)\times[-R^2,0]}\partial_tu_k^2dxdt&=&
      \mu_k^{\frac{4}{p-1}-(n-2)}\cdot\int^{t_k}_{t_k-\mu_k^2R^2}
      \int_{B_{\sqrt{\omega-t}}(x_k)\cap\Omega}\partial_tu^2(x,t)dxdt\\
      &=&o_k(1)
    \end{eqnarray*}
  when $1<p\leq p_s$. Consequently,
      $$
        \int^0_{-R^2}\int_{B_{b_\infty}(a_\infty)} v^2_t(x,t)dxdt=0
      $$
  for each $R>0$. Hence
      $$
        v(x,t)=v(x),
      $$
  i.e., $v(x,t)$ is a steady state satisfying \eqref{e0.18} in
  $Q^\infty$, which contradicts with \eqref{e0.19}. So, the conclusion is drawn. $\Box$\\

\vspace{40pt}

\section{Refined blowing-up III: Estimation on non-collapsing blowup set}

Let's start with the following criterion of blowup based on decay estimation in \cite{GK1} under critical and subcritical cases:

\begin{prop}\label{p4.1}
  Assume that $1<p\leq p_S$ and $\Omega$ is convex. Let $u$ be a maximal solution to \eqref{e1.1} on $\Omega\times[0,\omega)$. There exists a positive constant $\varepsilon_3$ depending only on $p$ and $N$, such that $(a,\omega)$ is not a blowing-up point, provided
    \begin{equation}\label{e4.1}
      r^{\frac{4}{p-1}-N}\int^{\omega-\frac{r^2}{e}}_{\omega-r^2}\int_{B_r(a)}|\nabla u|^2dxdt+r^{\frac{4}{p-1}-2-N}\int^{\omega-\frac{r^2}{e}}_{\omega-r^2}\int_{B_r(a)}u^2dxdt\leq\varepsilon_3
    \end{equation}
  holds for some small $r$ (say, $0<r\leq r_0$ with some $r_0$ depending on $p, N$ and $u_0$).
\end{prop}

\noindent\textbf{Proof.} Using case 2 in Theorem \ref{t3.1}, we want to show that \eqref{e3.25} holds for some function $\frac{1}{2}\leq\alpha(s)\leq1$.

 Setting
  $$
    I(r,\tau)\equiv\Big(\fint_{|y|=r}w(y,\tau)d\sigma\Big)^2\leq C_Nr^{-(N-1)}\int_{|y|=r}w^2(y,\tau)d\sigma,
  $$
 we have
   $$
    \int^1_{\frac{1}{2}}I(r,\tau)dr\leq C_N\int_{B_1\setminus B_{\frac{1}{2}}}\frac{w^2(y,\tau)}{|y|^{N-1}}dy\leq2^{N-1}C_N\int_{B_1}w^2(y,\tau)dy
   $$
 by Fubini's Theorem. Therefore, there exists $\frac{1}{2}\leq\alpha(\tau)\leq1$ such that
   \begin{equation}\label{e4.2}
     I(\alpha(\tau),\tau)\leq2^{N}C_N\int_{B_1}w^2(y,\tau)dy
   \end{equation}
 by intermediate value theorem for integral. Similarly, setting
   $$
     J(r,\tau)\equiv\Big(\fint_{|y|\leq r}wdy\Big)^2\leq C_Nr^{-N}\int_{B_r}w^2dy,
   $$
 we get
   \begin{equation}\label{e4.3}
     J(\alpha(\tau),\tau)\leq 2^NC_N\int_{B_1}w^2dy
   \end{equation}
 for $\frac{1}{2}\leq\alpha(\tau)\leq1$. Combining \eqref{e4.2} with \eqref{e4.3} and changing variables
  $$
   w\to(\omega-t)^{\frac{1}{p-1}}u,\ y\to\frac{x-a}{\sqrt{\omega-t}},\ \tau\to-\log(\omega-t),\ e^{-s}\to r^2,
  $$
\eqref{e3.25} holds true for our choice of $\alpha(s)$ using assumption \eqref{e4.1}. So, the conclusion is drawn. $\Box$\\

\noindent\textbf{Remark 4.1} The above proposition is in fact a special and slightly different version of Theorem 2 in \cite{CDZ} under critical and subcritical case. To estimate the blowup set for supercritical case, we need to quote it in full version with $p>p_S$:

\begin{prop}\label{p4.2}
   Assume that $p>p_S$ and let $u$ be a maximal solution to \eqref{e1.1} on $\Omega\times[0,\omega)$. There exist positive constants $\varepsilon_4, K>1, r_0<1$ depending only on $N, p, \Omega, \omega$, such that $(a,\omega)\in\overline{\Omega}\times\{\omega\}$ is not a blowup point, whenever
      $$
        r^{\frac{4}{p-1}-N}\int^{\omega-4r^2}_{\omega-9r^2}\int_{B_{Kr}(a)}(|\nabla u|^2+|u|^{p+1})dxdt<\varepsilon_4
      $$
   holds for some $0<r\leq r_0$.
\end{prop}

  The original version of Proposition \ref{p4.2} in \cite{CDZ} requires an assumption on convexity of $\Omega$. Utilizing the technics developed in \cite{CD}, it's not hard to see that a same result without convexity assumption also holds true. Now, let $u$ be a maximal solution of \eqref{e1.1} on $\Omega\times[0,\omega)$ and $\varepsilon_4$ be given in Proposition \ref{p4.2}, we define
   $$
     {\mathcal{S}}\equiv\Bigg\{a\in\overline{\Omega}\Big|\ r^{\frac{4}{p-1}-N}\int^{\omega-4r^2}_{\omega-9r^2}\int_{B_{Kr}(a)}(|\nabla u|^2+|u|^{p+1})dxdt\geq\varepsilon_4, \ \forall 0<r<r_0\Bigg\}
   $$
to be the blowup set of $u$ which is compact, and estimate its size for non-collapsing blowup:

\begin{theo}\label{t4.1}
  Assume that $p>p_S$. Let $u$ be a maximal solution blows up non-collapsing at $t=\omega$ and ${\mathcal{S}}$ be its compact blowup set whose Hausdorff dimension is no greater than $N-\frac{4}{p-1}-2$. We have $u$ is locally bounded near $(a,\omega)$ when $a$ doesn't belong to ${\mathcal{S}}\subset\overline{\Omega}$.
\end{theo}

\noindent\textbf{Proof.} For each $r_1<r_0$, noting first that ${\mathcal{S}}$ can be covered by closed balls $\overline{B_{Kr}(a)}$ with $a\in{\mathcal{S}}$. By Vitalli covering theorem, we can extract finitely many disjoint balls $\overline{B_{Kr_1}(a_k)}, k=1,2,\cdots, M$, such that
   $$
     {\mathcal{S}}\subset\bigcup_{k=1}^M\overline{B_{5Kr_1}(a_k)}
   $$
and
   $$
     r_1^{\frac{4}{p-1}-N}\int^{\omega-4r_1^2}_{\omega-9r_1^2}\int_{B_{Kr_1}(a_k)}(|\nabla u|^2+|u|^{p+1})dxdt\geq\varepsilon_4.
   $$
Therefore, taking any $q>1$, we have
   \begin{eqnarray*}
     \Sigma_{k=1}^M(5Kr_1)^{N-\frac{4}{p-1}-\frac{2q}{q-1}}&\leq&\varepsilon_4^{-1}(5K)^{N-\frac{4}{p-1}-\frac{2q}{q-1}}r_1^{-\frac{2q}{q-1}}\int^{\omega-4r_1^2}_{\omega-9r_1^2}\int_\Omega(|\nabla u|^2+|u|^{p+1})dxdt\\
      &\leq&\varepsilon_4^{-1}(5K)^{N-\frac{4}{p-1}-\frac{2q}{q-1}}\Bigg[\int^{\omega-4r_1^2}_{\omega-9r_1^2}\Bigg(\int_\Omega(|\nabla u|^2+|u|^{p+1}\Bigg)dx\Bigg)^qdt\Bigg]^{\frac{1}{q}}\\
      &\leq& C^*<\infty.
   \end{eqnarray*}
As the bound $C^*$ is independent of $r_1$ and $M$, letting $r_1\downarrow0$, we conclude that $\Big(N-\frac{4}{p-1}-\frac{2q}{q-1}\Big)$ dimensional Hausdorff measure is finite. So, the theorem has been proven by letting $q$ tends to infinity. $\Box$\\

It's interesting to ask whether the bound $N-2-\frac{4}{p-1}$ for Hausdorff dimension of the blowup set is optimal? We believe that it is true. However, even in critical case $p=p_S$, it is still open about the possibility of non-collapsing blowup comparing to that of finite time blowup for harmonic heat flow between surfaces \cite{CDY}.

\vspace{40pt}

\section{Examples of complete and incomplete blowups}

At the end of the paper, we present some new examples of complete and incomplete blowups.

\begin{theo}\label{t5.1}

 For $1<p<p_S$, there is no non-collapsing blowup, and hence all finite time blowup must be complete.

 For $p=p_S$ and convex domain $\Omega$, there is no type I non-collapsing blowup. Therefore, all finite time blowup must be complete in case $\Omega$ is a ball, $u_0=u_0(|x|)\geq0$.

 When $p>p_S$, the blowup must be non-collapsing under either one of the following conditions:

 (1)  $u$ is a borderline solution to \eqref{e1.1}, whose blowing-up in finite time has been proven in \cite{CD,CDZ} for star-shaped domain $\Omega$, or

 (2) $\Omega=B_R$, $u_0=u_0(r)>0$ and $u'_0(r)\leq0$ for all $0<r=|x|<R$.

On another hand, the blowup must be collapsing and hence complete under the following condition:

 (3) $\Omega=B_{R}\setminus B_{R_0}$, $u_0=u_0(r)$ for all $0<R_0<r=|x|<R$.\\
\end{theo}

\noindent\textbf{Proof.} At first, the conclusion for subcritical case has been proven in Corollary \ref{c2.1}. Secondly, when $p=p_S$, it was shown in \cite{MM} that all positive radial symmetric solution in a ball can only blow up in type I. Thus, a combination with Theorem 3.3 yields the result of complete blowup. Finally, for supercritical nonlinearity, one knows that all borderline solutions must blow up in finite time with non-collapsing energy in case $\Omega$ is star-shaped (see \cite{CDZ}\cite{CD}). So conclusion for (1) holds true. Another hand, under assumption (2), it was shown in \cite{FM} that
   $$
     \limsup_{t\to\omega^-}\int_\Omega|u|^{p+1}(t)<\infty.
   $$
Therefore, the blowup can not be collapsing in this case. At the end, we show that the blowup is collapsing under assumption (3). In fact, if it is not true, then the Hausdorff dimension of the blowup set ${\mathcal{S}}$ must be no greater than $N-2-\frac{4}{p-1}$. However, by radial symmetry, the Hausdorr dimension of ${\mathcal{S}}$ must be no less than $N-1$. Contradiction holds since
   $$
     N-2-\frac{4}{p-1}<N-1.
   $$
$\Box$\\

\vspace{40pt}

\section*{Acknowledgments}
The author would like to express his deepest gratitude to Professor Dong-Gao Deng, Xi-Ping Zhu, Kai-Seng Chou, Xu-Jia Wang and Neil Trudinger for warm-hearted helps. Special thanks
were also owed to his Ph.D advisor, Professor Kai-Seng Chou, for bringing
him to the topics and giving him many key observations.\\

\end{document}